# Tridiagonal and Block Tridiagonal computed Sparse Preconditioners for large Electrodynamic Electric Field Integral Equation (EFIE) Solution

Yoginder Kumar Negi, *Member IEEE,* N. Balakrishnan

*Abstract*— In this work, we propose simple and efficient tridiagonal computed sparse preconditioners for improving the condition number for large compressed Electric Field Integral Equation (EFIE) Method of Moment (MoM) matrix. The preconditioner computation is based on the triangle and block triangle interaction and filled tridiagonally. The computed preconditioner is highly sparse and retains the $O(N)$ complexity of computation and preconditioner matrix solution time. Numerical results show the efficiency and accuracy of the proposed preconditioner.

## I. Introduction

In the last few decades, Computational Electromagnetics (CEM) methods have gained popularity for various electromagnetic analyses due to their accuracy and efficiency. Frequency domain Electric Field Integral Equation (EFIE) based Method of Moment (MoM) [1] is one of the popular methods in CEM for solving complex electromagnetic radiation/scattering problems [2]. MoM leads to a dense matrix with $O(N^2)$ matrix fill time and memory requirement for $N \times N$ size matrix. Solving MoM system of equations requires $O(N^3)$ time with direct solver and $N_{itr}O(N^2)$ time with a conventional iterative solver for $N_{itr}$ iterations. The real world electromagnetics problems are geometrically large and complex; the solution of large scale problems with MoM is limited due to the high matrix storage, computation and solve cost. Direct solvers have the advantage of one-time factorization cost for a fixed time and memory. Solving large problems with a direct solver may be time-consuming and memory-intensive; even a few of the proposed fast direct solvers [3, 4] scale poorly for Thee-Dimension (3D) large complex problems. At the same time, the iterative solver needs less memory and fewer matrix operations than the direct solver. In the iterative solver, high storage and computation cost can be mitigated by incorporating matrix compression based fast solver methods like Multi Level Fast Multipole Algorithm (MLFMA) [5], pre-corrected FFT [6], Adaptive Cross Approximation (ACA) [7, 8], Hierarchal Matrices (H-Matrix)[9, 10] and IE-QR [11]. The matrix storage, fill-time and matrix-vector product time can be reduced to $O(NLogN)$ with a reduced solution time of $N_{itr}O(NlogN)$ for $N_{itr}$

This work was supported by the Department of Science and Technology (DST) Government of India under National Supercomputer Mission (NSM) project SP/DSTO-20-0130,

The authors Yoginder Kumar Negi and N Balakrishnan are with the Supercomputer Education and Research Center, Indian Institute of Science, Bangalore 560012, and India (email: yknegi@gmail.com, balki@iisc.ac.in).

iterations. EFIE being a Fredholm integral equation of the first kind, the eigenvalue tends to cluster at zero and infinity leading to the poor condition number of the matrix. As the number of unknowns increases, the number of large and small eigenvalues also increases, which leads to an increase in the ill-conditioning of the matrix. Ill-conditioned matrices are highly sensitive to perturbation in the system, which may jeopardize the accuracy of the solution and leads to a high iterative solution iteration count. Preconditioning [12, 13 and 14] of a matrix helps to improve the condition number of the matrix by clustering the eigenvalues around 1 and reducing the solution iteration count. Preconditioning is a way to convert the coefficient matrix from a system of the equation to the desired property system before the solution.

Broadly preconditioners can be classified as analytic and algebraic. Analytic preconditioners like Calderon preconditioner [15] are kernel dependent and sensitive to characteristics of the operator, thus applicable to a narrow class of problems. In comparison, algebraic preconditioners are more versatile and applicable for a broad range of problems. Incomplete LU (ILU) [16, 17, and 18] and Sparse Approximation Inverse (SPAI) [19, 20] are the few popular algebraic preconditioners for accelerating the iterative solution process. For significant size problems, ILU is limited due to the serial nature of LU factorization and selection of drop tolerance ($\tau$) and fill-in ($p$) parameters. In contrast, SPAI is limited by the quadric cost of computation and is applicable for parallel process-based matrix solutions. The Near-field matrix [21] of a fast solver can also be used as a preconditioner, but the high factorization cost limits the application as a preconditioner. A scaled near-field block-diagonal preconditioner is presented in [22, 23, 24, 25], but the diagonalization process is complex. The preconditioner should be simple and low cost in computation, and effective in improving the condition number of the matrix.

The diagonal and block-diagonal preconditioners are the simplest but are not effective in improving the condition number of the large size matrix [26]. In this paper, we propose novel sparse preconditioners based on the tri-diagonal and block tridiagonal interaction. The preconditioner is highly sparse and has a very low solution time. Tridiagonal matrix preconditioner is presented in [27, 28] and is applied for solving sparse matrices arising from Navier–Stokes equation. Our proposed sparse preconditioner is based on the triangle interaction and triangle cluster interaction at the lowest level of the binary-tree/oct-tree. The sparse preconditioners scale the columns of the coefficient matrix and improve the spectral property of the matrix, which further improves the iteration count during the solution process. The Numerical results show the accuracy and efficiency of the proposed preconditioner method.

The paper is organized as follows: in section II, a brief description of EFIE H-Matrix is presented. In section III, the proposed tridiagonal and block tridiagonal preconditioner is presented. In section IV, the efficiency and accuracy of the proposed H-Matrix are presented. Section V concludes the paper.

## II. EFIE H-MATRIX

EFIE based MoM is a popular method for solving open and closed conductor body problems in electromagnetics. For a 3D arbitrary shape conducting body, the EFIE boundary condition for an object illuminated with an incident field on the surface $S$ is given as:

$$|E_s(J) + E_i|_{tan} = 0 \quad (1)$$

where $E_s$ is the scattered electric field due to the induced surface current $J$ on the object illuminated by an incident electric field $E_i$ and $tan$ is Eq. (1), is a tangential component of the Electric Field. The scattered electric field can be further written as

$$E_s(J) = -j\omega A - \nabla \phi \quad (2)$$

where $A$ and $\phi$ represent the vector and scalar potentials, and $\omega$ is the angular frequency. Expanding scalar and vector potential and using the Galerkin testing method with RWG basis function [29], the resultant MoM system of Eq. (2)

$$[Z][x] = [b] \quad (3)$$

In Eq. (3), for a given unknown $N$, $[Z]$ is a dense MoM matrix of size $N \times N$, $b$ is an incident vector, and $x$ is a solution vector of size $N \times 1$. Dense MoM matrix leads to $O(N^2)$ matrix storage and filling time. The matrix storage and fill time can be reduced by incorporating fast algorithms for matrix filling and solutions. These methods work on the principle of analytic and algebraic matrix compressibility of far-field interaction blocks. MLFMA and FFT are analytic matrix compression methods. Algebraic matrix compression methods are ACA and H-Matrix, and IE-QR. These methods are kernel independent and easy to implement compared to an analytic method like MLFMA. In this work, we used half H-Matrix [30] with re-compressed ACA [31] to take advantage of algebraic compression and reduce the overall matrix storage requirement. For the H-Matrix construction, the compression scheme can be applied on a binary-tree based 3D geometry decomposition, where the matrix compression is applied for block interaction satisfying the admissibility condition.

$$\eta \, dis(\Omega_t, \Omega_s) \geq \min(dia(\Omega_t), dia(\Omega_s)) \quad (4)$$

The admissibility condition of Eq. (4) states that for matrix compression admissibility constant ($\eta$) times the distance between the test ($\Omega_t$) and source blocks ($\Omega_s$) should be greater than or equal to the minimum of the block diameter of the test block and source block. The binary-tree partition of the geometry is carried out until the number of elements in the block is less than or equal to 30 basis elements. At the leaf level, the block interaction not satisfying the admissibility condition is considered as a near-field interaction. In the case of the multilevel binary tree, the far-field block satisfying admissibility condition interacted at a higher level does not interact at the lower level. As the number of unknowns grows up, the solution time grows up due to the increase in iteration count. The condition number of the matrix deteriorates as the matrix size increases; furthermore, mesh inconsistencies and geometry type may also lead to matrix ill-conditioning. An ill-conditioned matrix leads to a high iteration count and solution time. Preconditioners can improve the condition number of the matrix and accelerate the solution time. In the next section, we propose a simple and efficient preconditioner based on the tridiagonal and block tridiagonal matrix computation.

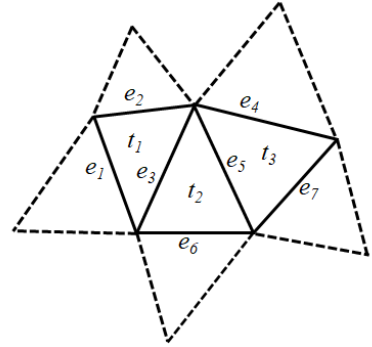

Fig 1. Triangle interaction for sparse matrix computation.

## III. TRIDIAGONAL AND BLOCK TRIDIAGONAL PRECONDITIONING

Iterative solution of large size matrix with Krylov subspace depends on the condition number of the matrix and matrix-vector product cost. Matrix-vector product cost can be reduced using fast solver methods, whereas iteration count in the matrix is condition number dependent and as the number of unknown grows, the condition number of matrix deteriorates. Preconditioning is one of the efficient methods to improve the condition number of matrix and expedite the solution process. The preconditioned system can be either as a left Eq. (5) or right Eq. (6), and is given as:

$$[P^{-1}][Z][x] = [P^{-1}][b] \quad (5)$$

$$[Z][P^{-1}][\tilde{x}] = [b] \quad (6)$$

Where $P$ is the preconditioner matrix, $Z$ is the EFIE MoM full matrix or compressed matrix, $b$ is excitation vector, $x$ and $\tilde{x}$ are the solution vectors, where $\tilde{x}$ is $[P]x$. To keep the cost of the iterative solution low, the preconditioner matrix should be highly sparse in nature and effective in improving the condition number of the matrix. In this section, we propose a new sparse preconditioner used as a left preconditioner for solving a large compressed matrix. The inverse of the MoM matrix is the ideal preconditioner for solving the MoM matrix with an iterative method, but the cost of the MoM inverse is memory and compute-intensive. Most of the algebraic preconditioners try to depict the inverse of the actual solution matrix. Therefore most of the proposed preconditioners in literature are derived from the actual matrix.

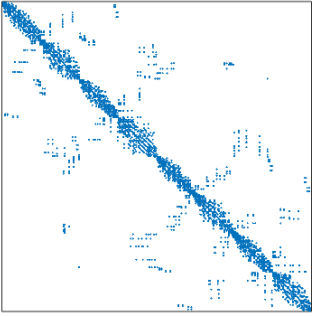

Fig 2. Sparse tridiagonal computed preconditioner for $1\lambda \times 1\lambda$ metallic plate.

The proposed first tridiagonally computed sparse preconditioner is highly sparse and is derived from the MoM matrix. The sparse preconditioner is computed by considering only the mesh triangle to triangle interaction in the MoM matrix. For an illustrative purpose, Fig. 1 shows 3 triangles that fill the MoM sparse preconditioner matrix. Here the triangle $t_2$ interacts with $t_1$ and $t_3$ to fill the 7 edges ($e_1, e_2 ... e_7$) indices in the MoM matrix. If the triangle $t_1$ is boundary triangle it interacts with $t_2$ else it as in the case of non-boundary, it interacts with $t_2$ and $t_0$. Dotted lines in the figure are the contributing triangles for the RWG edges. Fig. 2 shows the sparse tridiagonal preconditioner computed for $1\lambda \times 1\lambda$ metallic plate. The preconditioner is computed for 280 unknowns with 2688 Numbers of Non-Zeros (NNZ).

Similarly, we can compute block sparse preconditioner with triangle block interaction. Fast solvers like MLFMA and H-matrix relies on the oct-tree or binary-tree (Fig. 3) geometry division. For the divided geometry, the block matrix interaction is compressed at different levels with satisfying the far-field criteria. The non-far-field block interaction at the lowest level is considered near-field. Taking advantage of the geometric block partition for fast solvers, the preconditioner is computed for tridiagonal block interaction. The geometric block partition is done for triangles up to the desired level, with an average of 30 triangles in a group. For the block preconditioner computation as shown in Fig. 3, at the lowest level, triangles in block 1 interact with block 0 and 1, and in the case of boundary block 0, it interacts with block 1 triangles.

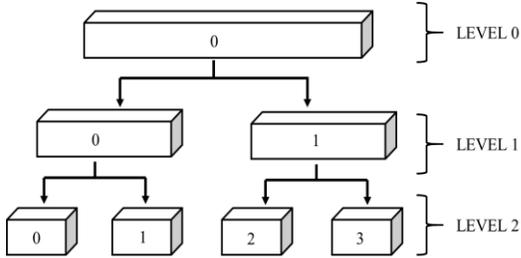

Fig. 3. Triangle block binary tree division of geometry for level 2.

Fig. 4 shows the sparse block tridiagonal preconditioner computed for $2\lambda \times 2\lambda$ metallic plate. The preconditioner is computed for 1160 unknowns at binary-tree level 5 with 156644 NNZ's.

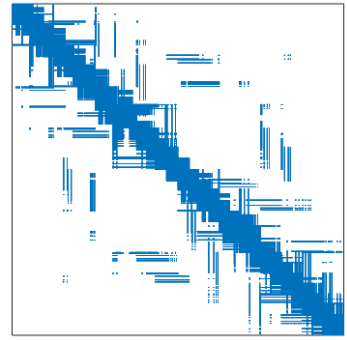

Fig. 4. Sparse block tridiagonal computed preconditioner for $2\lambda \times 2\lambda$ metallic plate for binary tree level 5.

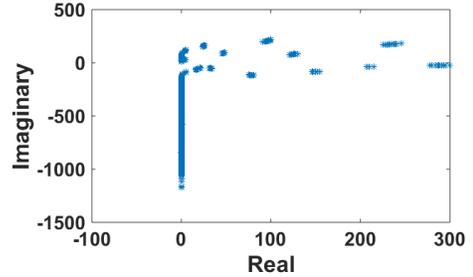

Fig. 5. Eigenvalue distribution of $1\lambda$ radius sphere MoM matrix with 5334 RWG edges

Fig. 5 above shows the eigenvalue distribution of $1\lambda$ radius sphere MoM matrix of size 5334×5334.

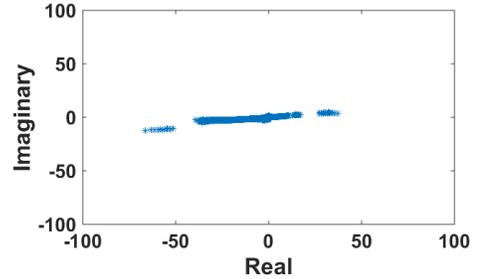

Fig. 6. Eigenvalue distribution after Tridiagonal computed sparse preconditioned $1\lambda$ radius sphere matrix

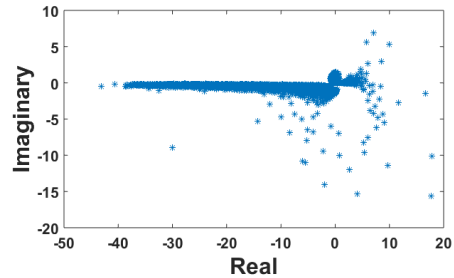

Fig. 7. Eigenvalue distribution after Block Tridiagonal computed sparse preconditioned $1\lambda$ radius sphere matrix

Fig. 6 and 7 shows the eigenvalue distribution of $1\lambda$ radius sphere with 5334 unknowns MoM matrix after preconditioning. It can be observed from the figure that the proposed preconditioners efficiently scales the columns of the MoM matrix and cluster the eigenvalues around 1, thus improving the spectral property of the EFIE matrix.

## IV. COMPLEXITY ANALYSIS

In this section, we show the efficiency of the proposed sparse preconditioners for set-up time; LU solve time and memory. The complexity analyses are carried for PEC plates with increasing unknown and size. One of the prime properties of a preconditioner should be its linear time complexity for set-up, and Fig. 8 shows that the proposed preconditioners retain the $O(N)$ complexity for computation.

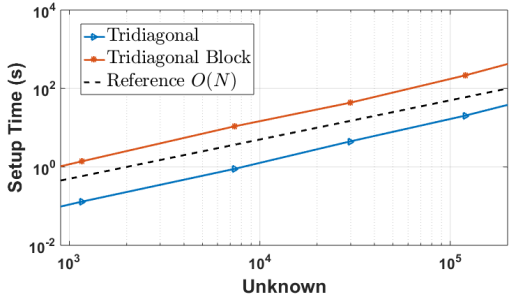
Fig. 8. Computation time for tridiagonal and block tridiagonal computed sparse preconditioner with increasing unknowns.

Iterative solver cost depends on the matrix-vector product time, and preconditioned iterative solver depends on the preconditioner LU solve time of the preconditioner. An efficient preconditioner should have a very less LU solve time with linear complexity. Fig. 9 below shows the linear LU solve time complexity for the proposed sparse preconditioners.

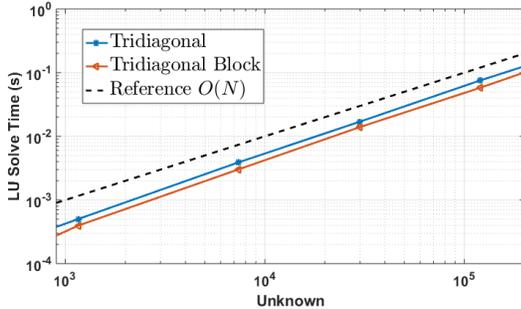
Fig. 9. LU solve time for tridiagonal and block tridiagonal computed sparse preconditioner with increasing unknowns.

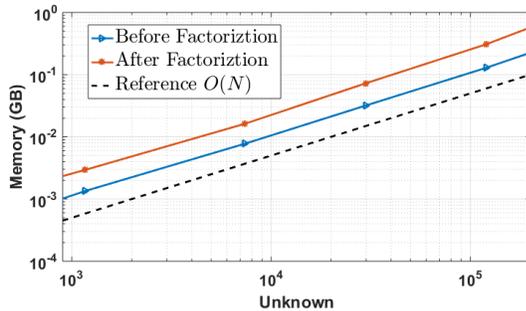
Fig. 10. Memory required for tridiagonal computed sparse preconditioner before and after factorization with increasing unknowns.

Along with time complexity, memory requirement plays a vital role in the efficiency of the solution process. A high memory preconditioner may jeopardize the iterative solution limiting the preconditioner applicable to a small size problem. Fig. 10 and 11 show the $O(N)$ memory complexity for tridiagonal and block tridiagonal computed sparse preconditioner before and after factorization.

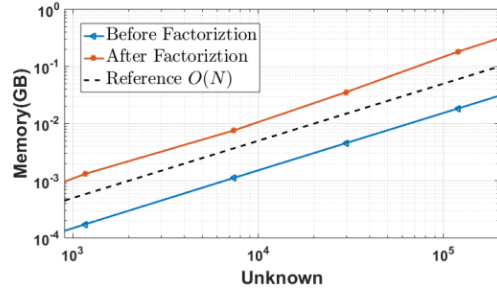
Fig. 11. Memory required for block tridiagonal computed sparse preconditioner before and after factorization with increasing unknowns.

## V. NUMERICAL RESULTS

In this section, we show the accuracy and the efficiency of the proposed preconditioners. All the simulations are done with ACA compressed H-matrix fast solver (compression tolerance = 1e-3) and solved with Krylov subspace-based iterative solver (GMRES) for convergence error of 1e-6 for Perfect Electric Conductor (PEC) geometry. Computation was carried out for double-precision data type on 128 GB memory and Intel (Xeon E5-2670) processor system.

### A. Accuracy

In the subsection, we show the accuracy of the proposed preconditioners for open and closed geometry. Fig 12 shows the monostatic RCS computation for $5\lambda$ square plate from MoM iterative solver, tridiagonal and block tridiagonal computed sparse preconditioned fast solvers for 7400 unknowns.

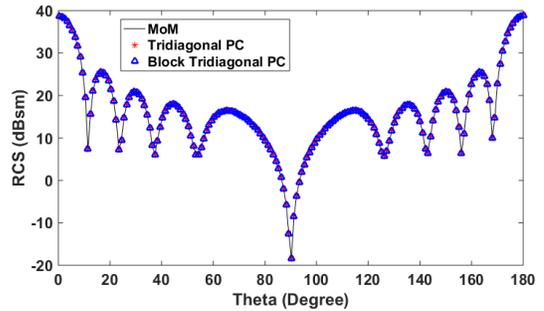
Fig. 12. Bistatic RCS of $5\lambda$ square plate with VV polarized plane wave incident at $\theta = 0°, \phi = 0°$ and observation angles $\theta = 0° \text{ to } 180°, \phi = 0°$.

TABLE I
PRECONDITIONER EFFICIENCY FOR DIFFERENT GEOMETRY

| | $N$ | PC Type | $T_{pc}$ (s) | $N_{itr}$ | $T_{pcsol}$ (s) | $T_{mmv}$ (s) | $T_{total}$ (H) (180 RHS) | Speed-up |
|---|---|---|---|---|---|---|---|---|
| 20λ Plate | 119,600 | TD | 63.8798 | 80 | 0.063675 | 0.703204 | 3.1384 | 1.46 |
| | | TD Block | 237.9965 | 51 | 0.067345 | 0.703204 | 2.5727 | 1.78 |
| | | ILUT | 273.7414 | 89 | 0.309779 | 0.703204 | 4.5838 | ---- |
| 5λ Sphere | 130,293 | TD | 75.9628 | 762 | 0.221020 | 1.450142 | 63.6923 | 1.36 |
| | | TD Block | 801.7663 | 711 | 0.079829 | 1.450142 | 54.6131 | 1.59 |
| | | ILUT | 842.6682 | 843 | 0.611366 | 1.450142 | 87.1266 | ---- |
| AC (1GHz) | 412,690 | TD | 724.8528 | 7788 | 0.520205 | 9.319170 | 3831.65 | 2.0 |
| | | TD Block | 6077.1354 | 7006 | 0.248953 | 9.319170 | 3353.40 | 2.3 |
| | | ILUT | 26138.339 | 12191 | 3.370137 | 9.319170 | 7742.02 | ---- |

The RCS is computed for the VV polarized plane wave incident at $\theta = 0°, \phi = 0°$ and observation angles $\theta = 0° \ to \ 180°, \phi = 0°$. It can be observed that the RCS computed from the preconditioned fast solvers agrees with the MoM computed RCS. For the MoM solution, the iterative solver takes 430 iterations, preconditioned tridiagonal preconditioner takes 45 iterations, and preconditioned block tridiagonal preconditioner takes 30 iterations to converge.

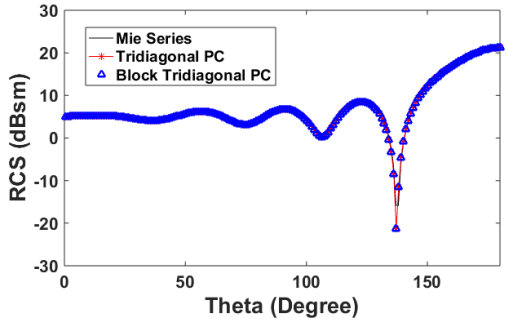

Fig. 13. Bistatic RCS of 1 λ radius sphere with VV polarized plane wave incident at $\theta = 0°, \phi = 0°$ and observation angles $\theta = 0°$ to $180°, \phi = 0°$.

Fig 13 shows the bistatic RCS computation for a 1 λ radius PEC sphere with 5334 unknowns. The RCS is computed using the Mie series analytic method and preconditioned tridiagonal and block tridiagonal fast solver for the VV polarized plane wave incident at $\theta = 0°, \phi = 0°$ and observation angles $\theta = 0° \ to \ 180°, \phi = 0°$. It can be observed that the RCS computed from the preconditioned fast solvers agrees with the Mie series RCS. For matrix solution, fast solver iterative solution takes 302 iterations, preconditioned tridiagonal computed sparse preconditioner takes 73 iterations and preconditioned block tridiagonal computed sparse preconditioner takes 52 iterations to converge.

*B. Efficiency*

In this subsection, the efficiency of the proposed preconditioners is validated. As discussed in our previous works [23, 25], the preconditioner efficiency cannot be concluded with fast iteration only. Along with iteration count, the preconditioner LU factorization solve time plays a vital role in overall solution time. In Table I, we show the solve time efficiency of our proposed preconditioners. The performance of the proposed preconditioners is compared with that of ILUT, with the parameters are chosen as given in [32]. The relative efficiency of a preconditioner depends on some key parameters: (a) $T_{pc}$: preconditioner set-up time, (c) $N_{itr}$: average number of iterations required for convergence for 1 RHS, (d) $N_{rhs}$: number of RHS vectors, (e) $T_{pcsol}$: preconditioner solve time and (f) $T_{mmv}$: MoM matrix-vector product time. The total solve time is given by:

$$T_{total} = T_{pc} + [N_{itr} \times N_{rhs} \times (T_{pcsol} + T_{mmv})] \quad (7)$$

Table I shows the speed-up efficiency of the proposed preconditioners for 180 RHS total solve time for a PEC plate, sphere and aircraft (AC). For plate and sphere, the iterations are computed for the VV polarized plane wave incident and observation angles at $\theta = 0° \ to \ 180° \ and \ \phi = 0°$ and for 14m length and 8m wingspan aircraft meshed with $\lambda/10$ element size at 1GHz the iterations are computed for the VV polarized plane wave incident and observation angles at $\theta = 90° \ and \ \phi = 0° \ to \ 180°$

## VI. CONCLUSION

The proposed preconditioners are simple to compute and effective in accelerating the iterative solution of large size problems. The preconditioners are sparse matrices based on the tridiagonal triangle and block tridiagonal interaction. These preconditioners maintain $O(N)$ set-up time, solution time and memory complexity. The preconditioners have a very low set-up time and can be divided into blocks and computed independently, thus making them highly efficient for parallel application.

ACKNOWLEDGMENT

We would like to thank the anonymous reviewers for the insightful comments and suggestions which helped us to improve the paper.

REFERENCES

[1] R. F. Harrington, *Field Computation by Moment Methods*, Krieger Publishing Co., Malabar, Florida, 1982.
[2] V. P. Padhy, Y. K. Negi and N. Balakrishnan, "RCS enhancement due to Bragg scattering," 2012 International Conference on Mathematical Methods in Electromagnetic Theory, pp. 443-446, 2012.
[3] Greengard, Leslie, Denis Gueyffier, Per-Gunnar Martinsson, and Vladimir Rokhlin. "Fast direct solvers for integral equations in complex three-dimensional domains." Acta Numerica, vol. 18, pp. 243-275, 2009.
[4] Shaeffer, John. "Direct solve of electrically large integral equations for problem sizes to 1 M unknowns." IEEE Transactions on Antennas and Propagation, vol. 56, no. 8, pp. 2306-2313, 2008.
[5] W. C. Chew, J. M. Jin, E. Michielssen and J. Song, *Fast Efficient Algorithms in Computational Electromagnetics*, Artech House, Boston, London, 2001.
[6] J. R. Phillips and J. K. White, "A precorrected-FFT method for electrostatic analysis of complicated 3-D structures," in *IEEE Transactions on Computer-Aided Design of Integrated Circuits and Systems*, vol. 16, no. 10, pp. 1059-1072, Oct. 1997
[7] M. Bebendorf, "Approximation of boundary element matrices," *Numerische Mathematik.*, vol. 86, no. 4, pp. 565–589, Jun. 2000.
[8] S. Kurz, O. Rain, and S. Rjasanow, "The adaptive cross-approximation technique for the 3-D boundary element method," *IEEE transactions on Magnetics*, vol. 38, no. 2, pp. 421–424, Mar. 2002.
[9] Hackbusch, Wolfgang. "A sparse matrix arithmetic based on H-matrices. Part I: Introduction to H-matrices." *Computing*, vol. 62, no. 2, pp. 89-108, 1999.
[10] W. Hackbusch and B. N. Khoromskij, "A sparse H-matrix arithmetic. Part II: Application to multi-dimensional problems", *Computing*, vol. 64, pp. 21-47, 2000.
[11] S. Kapur and D. E. Long, "IES3: Efficient electrostatic and electromagnetic solution*," IEEE Computer Science and Engineering* 5(4), pp. 60-67, Oct.-Dec. 1998.
[12] Benzi, Michele. "Preconditioning techniques for large linear systems: a survey." *Journal of Computational Physics,* Vol. 182, no. 2 pp. 418-477, 2002.
[13] Carpentieri, Bruno. "Preconditioning for Large-Scale Boundary Integral Equations in Electromagnetics." *IEEE Antennas and Propagation Magazine*, vol. 56, no. 6, pp. 338-345, 2002.
[14] Carpentieri, Bruno, Iain S. Duff, Luc Giraud, and M. Magolu monga Made. "Sparse symmetric preconditioners for dense linear systems in electromagnetism." *Numerical linear algebra with applications*, vol. 11, no. 8-9, pp. 753-771, 2004.
[15] Christiansen, Snorre H., and Jean-Claude Nédélec. "A preconditioner for the electric field integral equation based on Calderon formulas." *SIAM Journal on numerical analysis,* vol. 40, no. 3, pp. 1100-1135, 2002.
[16] Saad, Yousef. "ILUT: A dual threshold incomplete LU factorization." *Numerical linear algebra with applications* vol. 1, no. 4, pp. 387-402, 1994.
[17] Malas, Tahir, and Levent Gürel. "Incomplete LU preconditioning with the multilevel fast multipole algorithm for electromagnetic scattering." *SIAM Journal on Scientific Computing*, vol. 29, no. 4, pp. 1476-1494, 2007.
[18] Carpentieri, Bruno, and Matthias Bollhöfer. "Symmetric inverse-based multilevel ILU preconditioning for solving dense complex non-hermitian systems in electromagnetics." *Progress in Electromagnetics Research* vol. 128, pp. 55-74, 2012.
[19] Benzi, Michele, Carl D. Meyer, and Miroslav Tuma. "A sparse approximate inverse preconditioner for the conjugate gradient method." *SIAM Journal on Scientific Computing* 17, no. 5, pp. 1135-1149, 1996.
[20] Carpentieri, Bruno, Iain S. Duff, Luc Giraud, and Guillaume Sylvand. "Combining fast multipole techniques and an approximate inverse preconditioner for large electromagnetism calculations." *SIAM Journal on Scientific Computing*, vol. 27, no. 3, pp. 774-792, 2005.
[21] Carr, M., M. Bleszynski, and J. L. Volakis. "A near-field preconditioner and its performance in conjunction with the BiCGstab (ell) solver." *IEEE Antennas and Propagation Magazine,* vol. 46, no. 2, pp. 23-30, 2004.
[22] Y. K. Negi, N. Balakrishnan, S. M. Rao and D. Gope, "Null field preconditioner for fast 3D full-wave MoM package-board extraction," *2014 IEEE Electrical Design of Advanced Packaging & Systems Symposium (EDAPS)*, pp. 57-60, 2014.
[23] Y. K. Negi, N. Balakrishnan, Sadasiva M. Rao, and Dipanjan Gope. "Null-field preconditioner with selected far-field contribution for 3-D full-wave EFIE." *IEEE Transactions on Antennas and Propagation* 64, no. 11, pp. 4923-4928, 2016.
[24] Y. K. Negi, N. Balakrishnan, S. M. Rao and D. Gope, "Schur complement preconditioner for fast 3D full-wave MoM package-board extraction," *2016 IEEE Electrical Design of Advanced Packaging and Systems (EDAPS)*, pp. 163-165, 2016.
[25] Y. K. Negi, N. Balakrishnan, and Sadasiva M. Rao. "Symmetric near-field Schur's complement preconditioner for hierarchal electric field integral equation solver." *IET Microwaves, Antennas & Propagation* 14, no. 14, pp. 1846-1856, 2020.
[26] Gürel, Levent, Tahir Malas, and Özgür Ergül. "Preconditioning iterative MLFMA solutions of integral equations." In *2010 URSI International Symposium on Electromagnetic Theory*, pp. 810-813. IEEE, 2010
[27] Lukšan, Ladislav, and Jan Vlček. "Efficient tridiagonal preconditioner for the matrix-free truncated Newton method." *Applied Mathematics and Computation* 235, pp. 394-407, 2014.
[28] Li, Ji-cheng, and Yao-lin Jiang. "Generalized tridiagonal preconditioners for solving linear systems." *International Journal of Computer Mathematics* 87, no. 14, pp. 3297-3310, 2010.
[29] S. M. Rao, D. R. Wilton and A. W. Glisson, "Electromagnetic scattering by surfaces of arbitrary shape," *IEEE Transactions on Antennas and Propagation*, vol. 30, no. 3, pp. 409-418, May 1982.
[30] Y. K. Negi, "Memory Reduced Half Hierarchal Matrix (H-Matrix) for Electrodynamic Electric Field Integral Equation." *Progress In Electromagnetics Research Letters,* vol. 96, pp. 91-96, 2021
[31] Y. K. Negi, V. P. Padhy, N. Balakrishnan. "Re-Compressed H-Matrices for Fast Electric Field Integral Equation", *IEEE-International Conference on Computational Electromagnetics (ICCEM 2020)*, Singapore, August. 24-26, 2020.
[32] Li, X. Sherry. "SuperLU: Sparse direct solver and preconditioner.*" In 13th DOE ACTS Collection Workshop*. 2004.